\documentclass{article}
\usepackage[utf8]{inputenc}
\usepackage{amsmath}
\usepackage{graphicx}

\title{Winning Wordle Wisely\\
-or-\\
How to ruin a fun little Internet game with math}
\author{Martin B. Short}
\date{January 2022}

\begin{document}

\maketitle

\begin{abstract}
    In this missive, the author goes about ruining an online word guessing game called Wordle, by generally trying to mathify the whole thing.  On the other hand, one might consider that the author is trying to spice up mathematics by applying it to an online word guessing game.  Also, a spoiler warning: don't read past page 9 if you want to retain that Wordle magic.
\end{abstract}

\section{Introduction, and how to play Wordle}
It is now January 2022, and Wordle \cite{wordle} has taken the Internet by storm. The Google Trends timeline \cite{trends} currently shows Wordle at 100, after an apparently exponential rise in interest starting Jan 1, 2022.  No doubt there will be a similarly rapid drop in interest soon, given the general nature of Internet fads.\footnote{When I wrote this, I truly believed it.  However, it was recently announced that Wordle was purchased by the New York Times \cite{nyt} for an undisclosed (but seven figure!) sum.  So, I think it is likely that Wordle is here to stay in one form or another for the foreseeable future, even if drops out of the trending Internet conversation rather quickly.}  But, for now, the game is of great popular interest.

What is Wordle, you may ask?  It is a fun little Internet-based word guessing game.  The premise and gameplay of Wordle is simple.
    Any given game of Wordle is based around the player attempting to find a secret answer word through a sequence of guesses.
    The answer word is always a five letter word from the English language, at least according to some dictionary (more on this in a bit). 
    The player has a maximum of six guesses to try to find the secret answer word.
    Each of the player's guesses must also be a five letter word from the English language, again with the caveat about dictionaries from above.  I should also note that one aspect that makes Wordle especially fun is that every day there is only one unique Wordle for all players worldwide, which gives a certain social aspect to the rather simple game.

The way that the Wordle game provides information to the player so that they may make subsequently ``better'' guesses each time is by marking each letter in the most recent guess by one of three possible colors, each with a specific meaning:
\begin{enumerate}
    \item{Green (G):} a Green letter indicates that the secret answer word has that particular letter in that precise location.  This is a direct hit, if you will.
    \item{Yellow (Y):} a Yellow letter generally indicates that the secret answer word does have that letter somewhere, but not at that particular location.  This color is a little more complicated than Green, which is totally unambiguous, as there are edge cases to consider here. So, the full algorithm for yellow letters seems to be as follows.  First, remove any letters corresponding to Green letters from the secret answer, thus leaving a substring (this is so that letters already marked as Green won't also give Yellows).  Then going from left to right one guess letter at a time, check if that letter is located within the substring.  If so, then that guess letter will be marked as Yellow, and the letter within the substring that matched it is removed from the substring before checking the next letter in the guess (this is to prevent multiple Yellows coming off of a single letter within the answer).
    \item{Black (B):} Really Gray would be a better reflection of the color that Wordle uses here than Black, but I will use Black throughout this document.  Any letter not Green or Yellow is Black.  So, Black letters will generally indicate letters that are not contained anywhere within the secret answer.  The exception to this are letters that may appear, but for which there is already a corresponding guess letter labeled as Green or Yellow.  
\end{enumerate}

Here are a few examples of answers and guesses, and the corresponding letter color codes that Wordle would indicate, in case there is any confusion:
\begin{itemize}
    \item Answer: WEARY, Guess: TRAWL.  A straightforward example with resulting colors BYGYB.
    \item Answer: DRAFT, Guess: TWIST. Resulting colors are BBBBG.  Here, though the answer contains the letter T,  the first instance of T in the guess is given color B, since the T in the answer already corresponds to a letter colored G.
    \item Answer: GLEAN, Guess:  NANNY.  Resulting colors are YYBBB.  Only the first instance of letter N in the guess, read left to right, is colored Y.
    
\end{itemize}

In order to make the forthcoming discussion more meaningful, I should describe here a few decision and assumptions that will underlie everything I write going forward.
\begin{enumerate}
    \item Let's define the list of all words that Wordle will allow a player to guess to be the Wordle ``dictionary''.  This dictionary can be found online, as it turns out that one can simply inspect the Wordle source code and download it.  However, in much of my work in this paper, with exceptions specifically noted when they arise, I will not be using this ``true'' Wordle dictionary.  This is for a number of reasons.  First, the true dictionary can (and already did, after purchase by the NYT) change at any time, so the exact results one might obtain for a specific dictionary are not really that interesting; it is the methods and ideas that matter.  Second, the true Wordle dictionary is quite large, and contains many groupings of five letters that I simply refuse to acknowledge as English words, and which I would wager that 99.99\% of all native English speakers have never used, seen used in any context, or could even hazard a guess at the definition of.  So, I don't feel bad about excluding them, and for computational purposes I'd rather work with a smaller dictionary anyway.  And third, I would like to retain some mystery and challenge for myself while playing Wordle, and spending too much time analyzing the true dictionary might ruin this for me.  
Consequently most results I present here will be based off of a specific dictionary of five letter English words that I obtained using Mathematica 13.0.0.0.  After removing proper nouns from this dictionary and also removing a small number of words that contain non-standard English letters (letters with accents, etc), the dictionary contains $N=5170$ words (the ``true'' Wordle dictionary that I obtained on February 4, 2022 contains 12972 ``words'').  The variable $N$ shall be used henceforth to represent the size of the Wordle dictionary in use, whichever that may be.
    \item I will assume going forward that the secret answer to any given game of Wordle is drawn from the same full dictionary referred to above.  Further, I will assume that every word in the dictionary is equally likely to be the secret answer for any given game.  This is not true.  The secret answers to all Wordle games for the foreseeable future can be viewed through Wordle's code (but I would never dare to look at them myself), so they are clearly not being drawn randomly on each given day, and even if the sequence was originally randomly determined, it is doubtful that it contains any word more than once, so there is perhaps some memory to the random process.  Also, this list of foreseeable future answers is much smaller than the full dictionary, at only 2315 words (for the list I obtained on February 4, 2022, which I have not ever looked at, I swear), so it does not contain all dictionary words.  All this is ignored here.  It would be very easy to adapt any of the methods and algorithms I present here to use only this known set of possible answers, but if you are going to do that, why not just look at the source code for that day's solution and be done with it.
\end{enumerate}

Let us call the current list of words (as taken from the full dictionary or whatever subdictionary you cheaters might want to use) that are still potentially the secret answer for that game given the set of guesses made thus far and corresponding color codes received in response to them the ``viable set''.
Then each time a player makes a guess and receives in return a color code, the viable set can be updated, by either the game or the player.  If we denote the viable set at round $m$ (before the $m^\textrm{th}$ guess is made) by vector $\vec{v}^m$ (so each entry in this vector is a word from the dictionary), where at level $m=1$ the viable answer list is the full dictionary (by assumption), then the guess made and color code received at round $m$ will dictate what $\vec{v}^{m+1}$ is.  

Now, I have thought a fair bit about different ways of representing a game of Wordle, and how to most efficiently deduce the new viable set $\vec{v}^{m+1}$ from the prior viable set, the last guess, and color code received.  What I have settled on is, no surprise, a matrix.  Let us call this very helpful matrix the ``color code'' matrix, denoted by $C^m$, where the $m$ superscript again refers to the round of the game we are at before making guess $m$.  The columns of $C^m$ each correspond to a different possible answer word from the current viable set $\vec{v}^m$, while the rows of $C^m$ each correspond to all the unguessed words that can still be guessed from the full dictionary; the number of rows is $N+1-m$.  Matrix entry $C^m_{ij}$ is a representation of the color code that one would receive upon guessing word $i$ should the true answer be word $j$.  One could represent these color codes in a variety of ways, but for my purposes I will just use an integer from 1 to $M=3^5 - 5=238$, which encompasses all possible color codes (the 5 possible codes that are subtracted here represent codes where there are 4 Gs and 1 Y, which are not logically possible).  The ordering of the color codes among these integers generally does not matter, though I will make the decision that color code $M$ corresponds to GGGGG, for reasons that will be clear later. 

Then it is quite simple to update the viable set after making each guess.  If the guess made was word $i$, and the color code received was $c$, then the updated viable set is just the set of all columns in row $i$ of $C^m$ such that $C_{ij}=c$.
Of course, to use this method, one must precompute $C^1$ for a given dictionary; this need only be done once, and can then be stored for all subsequent uses.  On my (rather beefy) personal computer running all 12 cores in parallel with Matlab, it takes only 16 seconds to compute $C^1$ using the dictionary of size $N=5170$, and for the true Wordle dictionary of size $N=12972$ it takes a mere 100 seconds.  So, eminently doable.

 
 It is worth mentioning here that Wordle also has an optional ``hard mode'', in which each guess made must conform to all currently known information about the secret answer gleaned from prior guesses thus far.  That is, any guess made at level $m$ in hard mode must be a member of the current $\vec{v}^m$.  In order to adapt our method above to this mode, the only necessary change when going from $C^m$ to $C^{m+1}$ is to retain only those rows from $C^m$ that correspond to the words in $v^{m+1}$; in this way, $C$ will always be a square matrix.  We will discuss more about hard mode later.

So, having described the game, and how one might conceptualize it or implement it on a computer, we now come to the heart of the matter.  What is the ``best'' way to play Wordle?  In the subsequent sections I will provide some possible answers to this question, and raise several followup questions that would be interesting to tackle.  I should also note here that there are already many, many online sources (roughly 34.5M, according to Google \cite{gresults}) that discuss this question.  Given that the number of websites discussing this matter vastly outnumbers the size of the Wordle dictionary, probably what I write below is not unique.  However, unlike traditional academic research (which this article most certainly is not), I believe it is imperative that I \emph{not} read any of these articles/papers before presenting my own result into the world via this note.  This feeling I think can be attributed to a desire not to overly ``cheat'' at Wordle, which is a game that I quite enjoy after all.

\section{What is the point of Wordle?}
In general, the only goal of a given game of Wordle is to make sure that you find the secret answer word within your six allotted guesses.  So long as you do that, you ``win''.  Maybe, then, we should search for a method that guarantees that we will always win any game of Wordle, despite the secret answer word for that game?  This is a valid goal, but one I will not focus on too much.  The main reason being that it doesn't seem so hard to win at Wordle; not to toot my own horn, but I've never failed a game of Wordle myself (humble brag).  Of course, since Wordle is really using a subset of the full dictionary for its secret answers, and that subset is likely chosen to be words that are most familiar to people, this is not really a fair assessment of what might happen for me if the full dictionary were used as possible answers.  But, since the remainder of this work will be based on how an algorithm might best solve Wordle, familiarity with certain words or the lack thereof is not really relevant, so probably a computer can pretty easily win Wordle even if the full dictionary were used as possible answer.  In any case, this is something that I can and will assess for any algorithms described here.

I will take the position that the true goal of Wordle is not only to win, but to do so in the fewest guesses.  This is backed up by the game itself, which gives more impressive sounding commendations in response to wins when those wins take fewer attempts; from the lowly ``Phew'' on the sixth attempt, to ``Great'', then ``Splendid'', ``Impressive'', ``Magnificent'', to the ultra-elusive ``Genius'' if won on the opening guess (which is just 100\% pure luck, not Genius).  So, this is the goal I will focus on here.  If we achieve this goal, and find an algorithm  that, when it wins, takes very few guesses to do so, then there is a good chance that algorithm will always win (but this is not strictly guaranteed without further investigation).

Of course, the game of Wordle is also one of chance.  Unless you have managed to construct your guesses, along with the corresponding color codes, such that only one possible viable word remains (much more on this in a bit), your probability of selecting the correct word on your next guess is less than 1.  So, at best, we could state that the goal of Wordle is to find the solution to the game in the fewest \emph{expected} number of rounds, with expected here signifying expectation relative to some probability distribution.  This probability will of course take into account the fact that, in any given game, we don't know what the secret answer is, so must average over all the possible (remaining) answers.

Now that we have a goal in mind, we can elaborate on how one might achieve it \emph{in principle}. A complete strategy for Wordle would involve specifying, for any potential (achievable) color code matrix $C^m$, the next word one will choose $g^{m+1}$.  We will signify a given strategy by $S_k$, such that $g^{m+1}_k=S_k(C^m,m)$.
Of course the full space of all possible $C^m$ is enormous, and the full space of all possible strategies is bigger still.  For now, though, given the theoretical ability to delineate all possible strategies $S_k$, we could find the optimal such strategy in the following way.  For each such $S_k$, use that strategy against each possible hidden answer word $a_j$ for $j=1,2,\ldots,N$, in each case noting on which guess the answer is obtained: call this number $R_{kj}$.  In many cases the result will be $R_{kj}>6$, indicating that that particular strategy would fail against secret answer $a_j$, but this is ok, as strategies that display this property on any significant number of answer words will likely not be optimal by our measure anyway.  Then, for strategy $S_k$, the average number of guesses required across all possible answer words, again remembering that we assume all answer words are equally probable, is just 
\begin{equation}
    E_k = \frac{1}{N}\sum_{j=1}^N R_{kj}~.
\end{equation}
The optimal strategy(s) is then the strategy(s) $S_{k'}$ whereby $E_{k'}\leq E_k$ for all $k$; that is, the strategy with the smallest corresponding $E_k$.

Of course, while the above may be conceptually helpful in deciding what an optimal strategy would truly be, I cannot imagine that it is computable in this way, given the extremely large state spaces involved.  So, to try to find such an optimal strategy, we will need to construct some approximate algorithms and see if we can find optimal strategies using them instead.

\section{Approximate optimal strategies}
Here I outline a few methods of approximating the true optimality problem described above.
\subsection{Basic approaches}

\subsubsection{Most Rapid Decrease (MRD) algorithm}
Consider again that our goal is to solve Wordle in the least number of guesses, averaged over all possible secret answers.  Note also that, for every guess word that one makes in Wordle, the resulting color code obtained will result in a decreasing of the current viable set of possible answers, beginning from the full dictionary, as words get ruled out (columns get removed) by the information gained (some guesses/color codes may result in no reduction, but there is never an increase in the size of the viable set from a guess).  Then perhaps one approximate method of attaining our goal is to attempt the following: at every guess $m$, given the current color code matrix $C^m$, guess that word/row (or one of many if there are ties) that will result in the smallest size (in terms of the number of members) of the updated viable set $\vec{v}^{m+1}$ in expectation over all the current viable words from $\vec{v}^m$.

In this method, we are not really worrying about getting the answer correct with any of our guesses until we reach the point where the viable answer list has been pared down to 1 (or perhaps 2, at which point we simply guess one of them) remaining answers, and we are just aiming to achieve this as quickly as possible.  However, we could build into this method a simple sanity check: if, at any guess $m$, the current optimal set of guesses has some words overlapping with the viable set $\vec{v}^m$, then we will choose one of these overlapping words surely.  This is because these overlapping words not only accomplish the same viable set size decrease as the non-overlapping words, but also have some (maybe small) probability of being the correct answer themselves, thus possibly ending the game immediately upon their selection.

One thing to keep in mind with the MRD method is that it is a greedy algorithm \cite{greedy}.  This means that it does not look ahead beyond just the next step, and just tries to make the viable set as small as possible on the next step only.  This can clearly lead to suboptimality overall, as it is possible that a word that does not achieve MRD optimality when simply considering the next step could achieve a smaller expected viable set size when considering the next step and the step following it together.

Having described the MRD method, we now describe how it can be implemented computationally.  Luckily, it is extremely easy to do so, given our representation of the game by color matrix $C^m$.  First, note that for a given row $i$ of $C^m$, the probability that we will receive a given color code $c$ upon guessing word $i$ is found by simply counting the number of times $c$ appears in that row, then normalizing by the total number of columns.  This is because we assume each word in the viable set is equally probable to be the secret answer.  Let us denote the number of times that color code $c$ appears in row $i$ of $C^m$ by $L(c,i,C^m)$, and denote the number of columns of $C^m$, which is the size of the current viable set, by $N_v$.  But also note that the size of the updated viable set $v^{k+1}$, upon guessing word $i$ and receiving color code $c$, is also $L(c,i,C^m)$.  Then if we guess word $i$, the expected length of the viable set $v^{k+1}$, which we will denote as $\overline{L}_i$ is simply given by
    \[
    \overline{L}_i=\frac{1}{N_v-1_{M\in C_i^m}}\sum_{c=1}^{M-1}L^2(c,i,C^m)~,
    \]
where $1_{M\in C_i^m}$ is an indicator function that tells us whether or not row $i$ contains color code $M$, which is GGGGG, which can by definition only appear at most once in a row, and only for those words that are in the viable set.  We also therefore excluded the color code GGGGG in the sum.  The reason why I have chosen to remove GGGGG in this case is that, if the color code we receive is GGGGG, then the game has ended and we have already won, so computing $\overline(L)$ is nonsensical.  Plus, our sanity check that will prefer optimal words that are members of the viable set over those that are not will cover this case explicitly in order to break ties.

After computing the $L_i$ for all rows of $C^m$, the MRD optimal next guess is then taken from that list of words $g_{i'}$ such that $\overline{L}_{i'}\leq \overline{L}_i$ for all $i$.  If there is any overlap between this MRD optimal list and $\vec{v}^m$, then that overlapping subset is the MRD optimal list.

Because all that is required to find the MRD optimal guess is to count how many times each color code appears within each row of the matrix $C^m$, then square each count and sum them all up, it can be done extremely quickly on a computer.  I have found that this algorithm can be applied in real time on my personal computer to determine the MRD optimal next guess at any point in a game of Wordle almost instantly.

Given the MRD algorithm, one might naturally wonder: what is the best \emph{opening} guess for Wordle?  This is the only guess that does not vary from one game to the next, as the viable list $\vec{v}^1$ is always the full dictionary.  Well, given the full dictionary I am using as described above, the MRD optimal opening guess word is... drumroll please...

\newpage
\begin{center}
    TARES
\end{center}

\noindent Yes, TARES.  As in, ``Every morning he TARES his scale when making his espresso.'' 

Now, if you instead use the true Wordle dictionary to answer this same question, one obtains the word LARES.  Yes, LARES.  As in, ``I have no idea what LARES means, but it only differs from TARES by one letter.''  Note that LARES is not present in the smaller dictionary.  Also, TARES is the third best word by this measure in the true Wordle dictionary (number two is RALES).

\subsubsection{Greatest Expected Probability (GEP) algorithm}
One drawback of the MRD method is that it generally foregoes the possibility of randomly selecting the answer word in each round by drawing its optimum from the full dictionary (previously guessed words excluded).  While this may not give up much if the viable possible answer list still contains tens of words at a given point, it can make a bigger difference later on when the viable word lists become much smaller.  It is also not directly applicable when playing on hard mode.  One could adapt MRD to hard mode by simply adjusting $C^m$ for hard mode as described above.  However, given that each guessed word on hard mode has the potential to be the correct answer, one might consider if another approximate optimality algorithm might be more appropriate than MRD in this case.

Here, we consider such an alternative, which I will call the Greatest Expected Probability (GEP) algorithm.  In this approximate problem, the goal is not to pick that guess that will lead to the shortest expected updated viable set, but to instead pick that word that leads to the greatest expected probability that you could randomly choose the correct answer on the next round of play.  There is of course a connection between these probabilities and the lengths of the potential updated viable sets $L(c,i,C^m)$ discussed above.  Specifically, let us now define  $P(c,i,C^m)=1/L(c,i,C^m)$, where the interpretation of this quantity is: if we choose for our guess word $i$ and receive color code $c$ (and assuming it is not GGGGG), then the probability that our next guess will be the correct answer (assuming we will only choose guesses from the updated viable answer list) is $P(c,i,C^m)$.  Of course, this definition is only valid for color codes $c$ that actually appear within row $i$ of $C^m$. Then we can simply find for each guess word $i$ the expected such probability via
\begin{align}
\overline{P}_i&=\frac{1}{N_v-1_{M\in C^m_i}}\sum_{c=1,c\in C^m_i}^{M-1} L(c,i,C^m)P(c,i,C^m)~\textrm{or}\\
\overline{P}_i&=M_i(C^m)/(N_v-1_{M\in C^m_i}),
\end{align}
where $M_i(C^m)$ is the number of unique color codes (not counting GGGGG) that appear in row $i$ of matrix $C^m$.  The GEP optimal guess is then a member of that set of guesses $g_{i'}$ such that $\overline{P}_{i'}\geq \overline{P}_i$ for all $i$. As with MRD above, if there are any GEP optimal words that are also members of the current viable set, we will certainly choose one of those.

Now, we can ask what word is the GEP optimal opening guess for Wordle.  The answer is...

\newpage
\begin{center}
    TARES
\end{center}

\noindent Yes, TARES.  As in, ``The GEP optimal opening word for Wordle is TARES, just like in MRD!'' Here, I will not qualify this answer by noting which dictionary this pertains to, as it is the GEP optimal in both!

Of course, it is not generally the case that $\overline{P}_i = 1/\overline{L}_i$, so the minimum $\overline{L}$ did not have to correspond to the maximum $\overline{P}$, and indeed does not in the true Wordle dictionary.  To get a feel for the overlap between the best words for the two algorithms, I list the top ten best words (See Table \ref{tab:topten}) for both MRD and GEP (for the smaller dictionary) below, and there is no overlap outside of the top spot, with the exception that TEARS appears in both lists.
\begin{table}[h!]
    \centering
    \begin{tabular}{r|c|c|c|c}
         Rank & MRD & $\overline{L}_i$ & GEP & $\overline{P}_i$  \\
         \hline
         1 & TARES & 117.54 & TARES & 0.0354\\
         2 & RATES & 120.07 & TEARS & 0.0350\\
         3 & TALES & 122.34 & TIRES & 0.0344\\
         4 & ALOES & 122.37 & TRIES & 0.0340\\
         5 & SANER & 125.87 & PARES & 0.0340\\
         6 & ROLES & 128.50 & TALES & 0.0337\\
         7 & LANES & 131.13 & CARES & 0.0337\\
         8 & RILES & 133.90 & PEARS & 0.0337\\
         9 & TEARS & 135.12 & PORES & 0.0337\\
         10 & ROTES & 136.32 & SLATE & 0.0335\\
         5169 & FUZZY & 2148.3 & JAZZY & 0.0074\\
         5170 & YUKKY & 2161.1 & PZAZZ & 0.0056\\
         
    \end{tabular}
    \caption{The top 10 best Wordle opening words, according to MRD and GEP optimality.  For fun, I have also included the two \emph{worst} opening words in each case.  No, I wasn't aware that PZAZZ was a word either.}
    \label{tab:topten}
\end{table}

\subsubsection{Generalizing from there}
While the MRD and GEP algorithms were motivated by specific, interpretable considerations about what we might want in a Wordle guess, such limitations have never stopped a mathematician before.  Hence, we will now consider a more general optimal algorithm that encompasses both MRD and GEP, and which I will refer to as $p$-optimality, defined by the quantity
\begin{equation} \label{eq:popt}
\overline{f}(p)_i=\frac{1}{N_v-1_{M\in C^m_i}}\sum_{c=1,c\in C^m_i}^{M-1} L(c,i,C^m)L^p(c,i,C^m)~.
\end{equation}
Here, the eponymous parameter $p$ could in general be any real number we like\footnote{As an interesting mathematical aside, for $p\geq 0$ these $\overline{f}(p)_i$ are directly related to the vector $(p+1)$-norms of the potential viable set length vectors.}.  Since in general we would probably prefer guesses that lead to new viable sets that are small in size, if $p>0$ we will choose the $p$-optimal guess to be that word $i'$ that minimizes $\overline{f}(p)_i$ over all $i$, while if $p<0$ we will choose the word that maximizes $\overline{f}(p)_i$.  As always, we will always check if any of the words in this optimal set is a member of the current viable set, and prefer those if so.  Then it is clear that if we choose $p=1$ we recover MRD and if we choose $p=-1$ we recover GEP.

Given the general nature of the $p$-optimality method, the sky's the limit, and we can try out many different $p$ values and see which is best!  Of course, for each $p$ there may be a different optimal opening Wordle guess.  But there are some limits we can examine that will tell us what to expect here.  First, consider the case where $p\to -\infty$.  Here, the only terms that will materially contribute to the sum in \eqref{eq:popt} are those for which $L(c,i,C^m)=1$.  Since we are trying to maximize $\overline{f}$ in this case, this particular $p$-optimality is just trying to find that word that has the most potential viable sets of size one that arise from it.  In the smaller dictionary this word is PLATS.  Yes, PLATS.  As in, ``I think I know how to use PLATS in a sentence, but now that I'm on the spot, I find I can't confidently do so.''  On the other end of the spectrum is $p\to \infty$, in which case by far the most important term in each sum is that entry with the largest value of $L(c,i,C^m)$.  In this case we are trying to minimize $\overline{f}$, so this particular $p$-optimality is just trying to find that word whose biggest possible resulting viable set is the smallest; basically minimizing the worst-case-scenario damage.  In the smaller dictionary this word is ALOES.  Yes, ALOES.  As in, ``the plural of aloe is ALOES.''

Now, what would be really spiffy would be to test out different values of $p$ and see which one leads to the best overall outcomes in actual games of Wordle.

\subsubsection{Testing out the algorithm(s)}
In order to give a full assessment of these algorithm(s) as potential optimal algorithms for Wordle, I have tested each method for every possible game of Wordle that could occur, in each case using the appropriate optimal opening word.  That is, for some specific values of $p$, including of course -1 and 1 to include GEP and MRD, I have first determined the optimal opening word and then have simulated games of Wordle (both in normal and hard modes) for all possible secret answers, using $p$-optimality to select the ``best'' guess at each round in each game, then measured how many rounds $R_{pj}$ are needed to win for each secret word $j$.  This might sound like it would take a long time to do, but it only takes around 30 seconds for each $p$ value (and less for hard mode).

Before divulging the results, I should offer a few finer points of clarification.  First, I should describe how I handle situations where there are multiple optimal guesses in a round.  For example, suppose I am using $p=0.25$ and there are two words that both have the same minimal $\overline{f}(0.25)$, and I must choose one of them for my next guess.  As a first step, since it is quite easy to compute the various $\overline{f}(p)_i$ for many different $p$ values simultaneously, I will always compute using at least $p=-1$, $1$, $-10$, and $10$, and will break any tie for the actual $p$ value in question using the results for these other $p$ values in the order I have presented them here.  So, in this specific example of a tie for $\overline{f}(0.25)$, I would choose among those tied words the one(s) that have the optimal $\overline{f}(-1)$, then break any remaining tie among those words using the optimal $\overline{f}(1)$, then $\overline{f}(-10)$, then $\overline{f}(10)$.  Admittedly these specific four $p$ values (and their ordering here) are somewhat arbitrary, but they were all chosen to represent easily interpretable optimal properties, as discussed above (since these options include MRD and GEP, and $\pm 10\approx \pm\infty$ as far as I'm concerned).  In the case that all of this still does not settle the tie, then I deviate from $p$-optimality and instead just choose the word that has the highest expected quantity of G and Y colors in the resulting viable set.  This choice is partly because I had to choose something, and partly because these are the types of answers I like to see when playing Wordle (more on that later).  If that still doesn't work, then I try just for the highest expected number of G colors, for similar reasons.  Finally, if all that fails to break the tie, I simply choose the first word alphabetically; this does happen, and is pretty frequent in hard mode (much less so in normal mode), since all the possible choices are by definition pretty similar.

Second, I have chosen to iterate guessing for each of these simulated games until one of three things happens.  The first is that I happen to guess the correct word, in which case I of course stop and note the round I just won on.  The second is that I narrow the viable set down to only a single possible word, in which case the optimal answer is obvious and I can stop iterating, but still note that it will take me one more round to win.  The final case is that I have narrowed the viable set down to two words.  Here, it is clear that the best way to proceed is to select one of the two words at random (as neither can be better than the other), in which case I either win immediately or I narrow the viable set down to one word and win in one additional round.  The upshot is that, in expectation after achieving a viable set of size two, I will win in 1.5 more rounds.  Along these lines, I did not restrict my simulations to a maximum of six rounds in order to win, as I wanted to track how many rounds each would take even if they had unlimited tries.  Of course, any number of rounds $R_j>6.5$ is a definite loss for the algorithm, while a value of $R_j=6.5$ indicates that after the fifth guess it had the options narrowed down to a set of two, so it was still possible to win at least.
\begin{table}[ht!]
    \centering
    \begin{tabular}{r|c|c|c|c|c|c}
         $p~~$ & Opener & Mean $R^n_{pj}$ & $R^n_{pj}=6.5$ & $R^n_{pj}=7$ & Mean $R^h_{pj}$ & $R^h_{pj}\geq 6.5$  \\
         \hline
         -2.00 & PLATS & 3.8097 & 6 & 3 & 3.9162  & 207  \\
         -1.75 & PLATS & 3.8070 & 4 & 3 & 3.9137 & 208\\
         -1.50 & PLATS & 3.8039 & 4 & 3 & 3.9083 & 203\\
         -1.25 & TEARS & 3.8236 & 16 & 8* & 3.9431 & 249\\
         -1.00 & TARES & 3.7752 & 6 & 2 & 3.8781 & 183\\
         -0.75 & TARES & 3.7712 & 6 & 2 & 3.8737 & 183\\
         -0.50 & TARES & 3.7739 & 2 & 2 & 3.8760 & 189\\
         -0.25 & TARES & 3.7733 & 4 & 2 & 3.8731 & 190\\
         0.25 & TARES & 3.7843 & 4 & 0 & 3.8903 & 200\\
         0.50 & TARES & 3.7901 & 4 & 0 & 3.8959 & 196\\
         0.75 & TARES & 3.7983 & 2 & 0 & 3.9056 & 200\\
         1.00 & TARES & 3.8008 & 2 & 0 & 3.9124 & 208\\
         1.25 & TARES & 3.8056 & 2 & 0 & 3.9186 & 208\\
         1.50 & ALOES & 3.8851 & 2 & 3 & 4.0451 & 286\\
         1.75 & ALOES & 3.8905 & 2 & 3 & 4.0576 & 306\\
         2.00 & ALOES & 3.8940 & 2 & 3 & 4.0578 & 305\\
    \end{tabular}
    \caption{The performance of the $p$-optimality algorithm for various $p$ values, in both normal mode (denoted by superscript $n$) and hard mode (denoted by superscript $h$).  The columns labeled $R_{pj}=X$ give the frequency of observing that particular $R_{pj}$.  Note the * for $p=-1.25$ in the $R^n_{pj}=7$ column: this denotes that there were actually 6 instances of $R^n_{pj}=7$ and 2 instances of $R^n_{pj}=7.5$ (all other rows had a maximum of $R_{pj}^n=7$). }
    \label{tab:pperform}
\end{table}

I present results comparing performance across several different $p$ values in Table~\ref{tab:pperform}.
The table presents the different $p$ values in two lights: one is the average number of rounds it takes to win and the other is a measure of how often the algorithm fails (or potentially fails) to win at all.  One particular $p$ value that stands out as quite terrible in both regards is $p=-1.25$, which has some of the most (and worst) losses, and also has a strangely high number of expected rounds, at least in comparison to its neighboring $p$ values.  Of note here is that this value is also the only one that uses an opening word that seems to come out of nowhere: TEARS.  Yes, TEARS.  As in, ``the performance of the algorithm at $p=-1.25$ brings TEARS to my eyes.''  All other values use the magical TARES or use PLATS (for $p\to -\infty$ as discussed above) or ALOES (for $p\to\infty$ as discussed above).  So, I'm blaming the strange and terrible performance here on the opening word.

\begin{figure}
\includegraphics[width=0.48\textwidth]{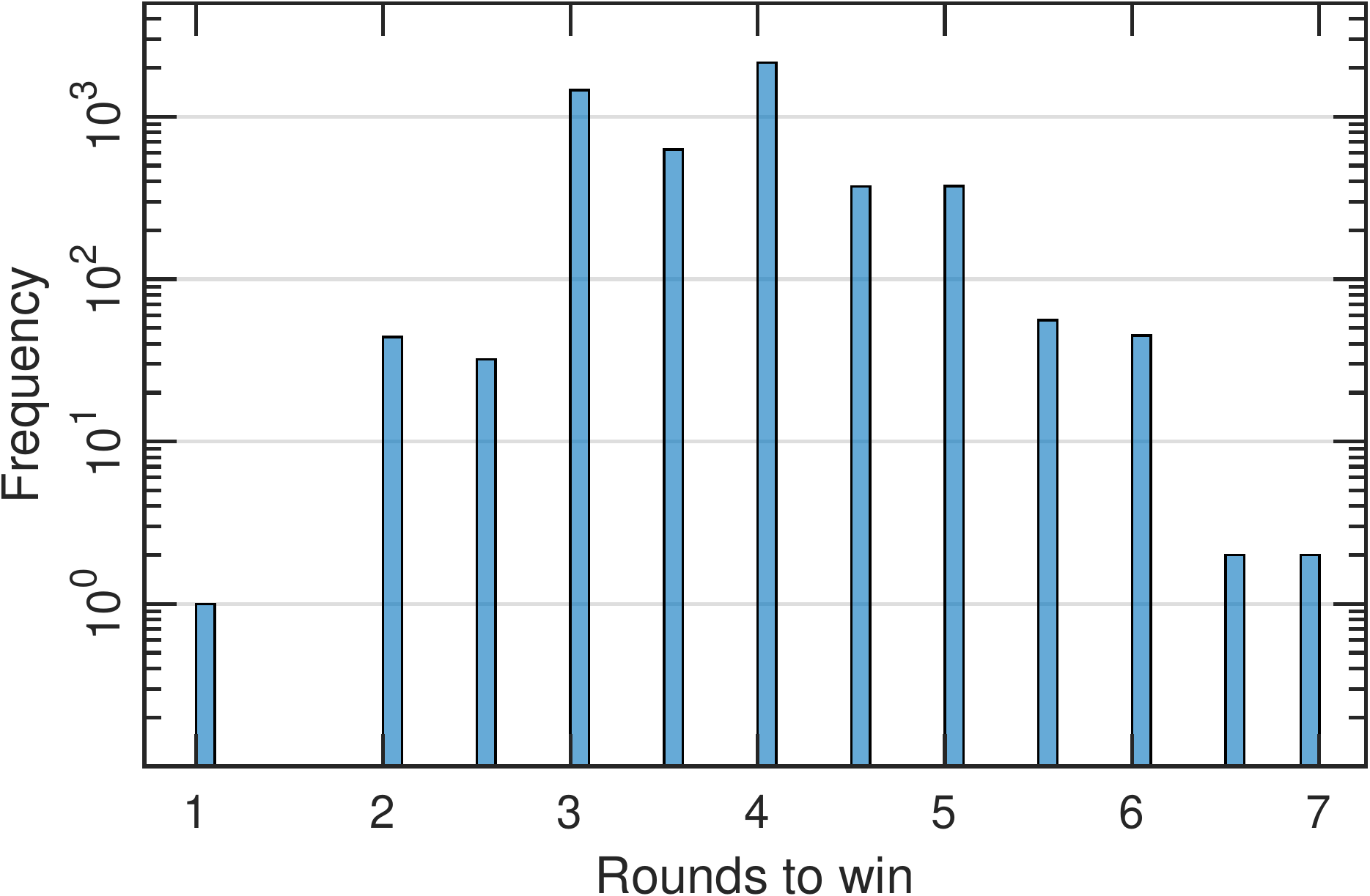}\hfill
\includegraphics[width=0.48\textwidth]{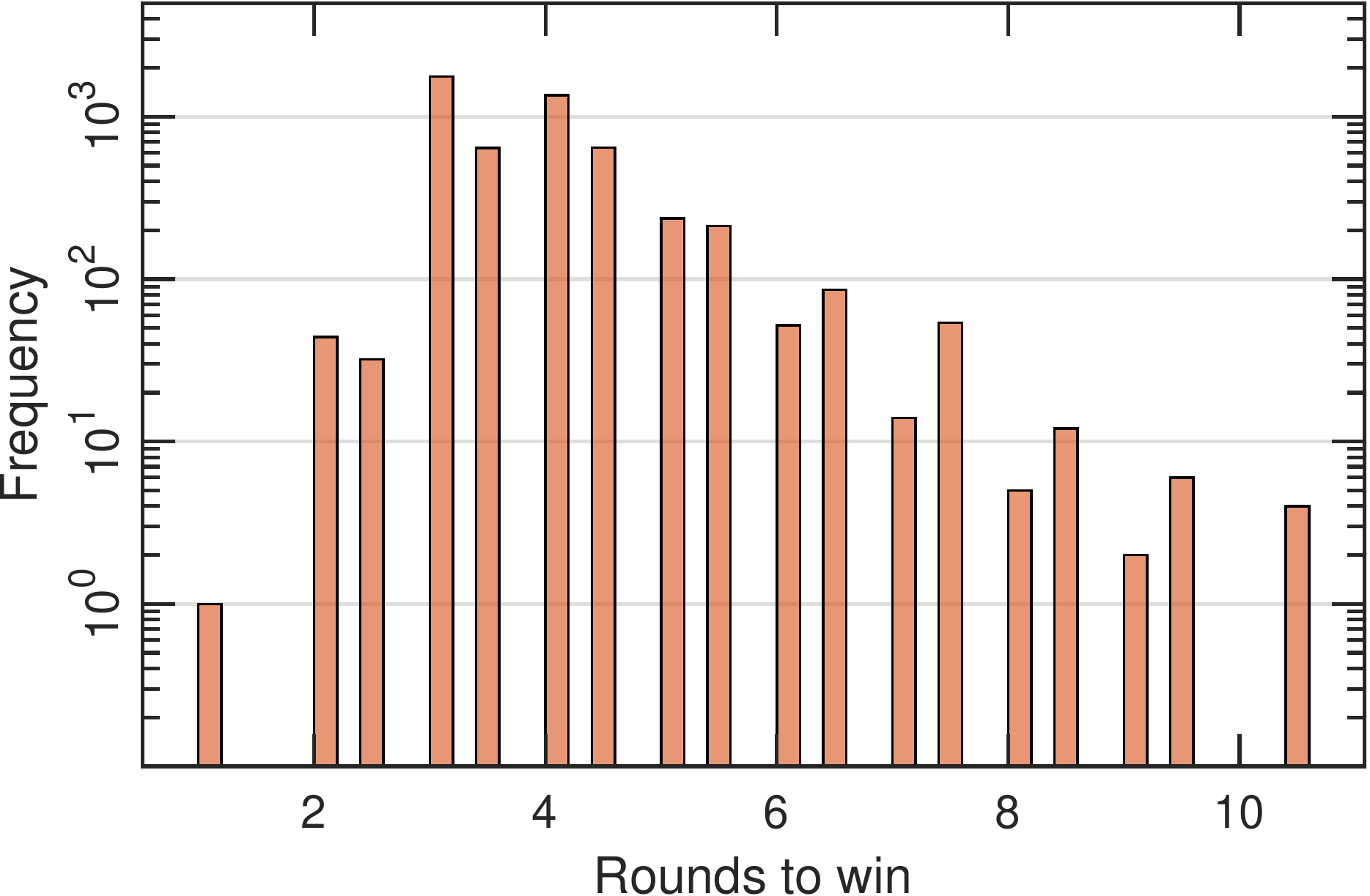}
\caption{The number of rounds required to win every possible game of Wordle using the smaller dictionary in (left) normal mode using $p=-0.5$ optimality and (right) hard mode using $p=-0.75$ optimality, in each case beginning from the optimal opening word TARES. \label{fig:phist}}
\end{figure}

Let us now focus our discussion on normal mode.  An interesting finding here is that, for intermediate $p$ values, there seems to be a tradeoff happening between losses and average rounds needed.  For example, the lowest average rounds was for $p=-0.75$, but it also definitely lost twice and possibly lost six times.  On the other hand, $p=0.75$ never lost for sure, and only possibly lost two times, but had a somewhat higher average number of rounds. A comparison between GEP and MRD also reveals this tradeoff, with MRD taking more rounds to win but losing less frequently.  The value $p=-0.5$ seems to offer a nice middle ground here, with a nearly minimal average number of rounds and only losing at most four times.  To provide a little more detail, I have singled out the specific value $p=-0.5$ as a strong performer and plotted a histogram of $R_{pj}$ in this case; see Fig.~\ref{fig:phist}.

In case you're interested, and why wouldn't you be, there is a clear trend in terms of the words that lead to fail states for the algorithm.  There are 40 distinct words that end in failure, with the most common failed word across all $p$ values being FAZES, followed by FAXES and HAZES.  Generally speaking, many (14/40) of these fail words follow the pattern \_A\_ES, including those words that were the only two losers in the cases with the least losses (BABES and HAZES).

Hard mode presents a somewhat different picture.  Here, $p=-0.75$ is almost an unambiguous winner, as it ties for least losses (with GEP) and has nearly the smallest expected number of rounds, beaten ever so slightly by $p=-0.25$.  A histogram of results for $p=-0.75$ is shown in Fig.~\ref{fig:phist}.  Unsurprisingly perhaps, GEP definitely beats MRD in hard mode; GEP was designed specifically with hard mode in mind, after all.  Other than that, there isn't too much to say about hard mode other than it is, in fact, hard, relative to normal mode at least.

To give a nod to the true Wordle dictionary, I have calculated at least the normal mode results equivalent to those in Table \ref{tab:pperform} for this larger dictionary.  Results are shown in Table \ref{tab:fullperform}.
\begin{table}[ht!]
    \centering
    \begin{tabular}{r|c|c|c|c}
         $p~~$ & Opener & Mean $R^n_{pj}$ & $R^n_{pj}=6.5$ & $R^n_{pj}\geq 7$ \\
         \hline
         -2.00 & PLEAT & 4.1856 & 62 & 58 \\
         -1.75 & PLEAT & 4.1691 & 62 & 39 \\
         -1.50 & PEATS & 4.1600 & 50 & 39 \\
         -1.25 & PELAS & 4.1464 & 48 & 61 \\
         -1.00 & TARES & 4.0889 & 36 & 16 \\
         -0.75 & TARES & 4.0859 & 36 & 15 \\
         -0.50 & TARES & 4.0807 & 26 & 10 \\
         -0.25 & TARES & 4.0827 & 30 & 13 \\
         0.25 & TARES & 4.0910 & 42 & 6 \\
         0.50 & LARES &   4.1156 & 22 & 13 \\
         0.75 & LARES &   4.1240 &24 &  13\\
         1.00 & LARES & 4.1263 & 26 & 10 \\
         1.25 & LARES & 4.1327 & 24 & 10\\
         1.50 & LARES & 4.1378 & 26 & 12 \\
         1.75 & LARES & 4.1452 & 28 & 12 \\
         2.00 & LARES & 4.1479 & 28 & 12 \\
    \end{tabular}
    \caption{The performance of the $p$-optimality algorithm for various $p$ values in normal mode (denoted by superscript $n$) when used on the true Wordle dictionary.  The columns labeled $R_{pj}=X$ give the frequency of observing that particular $R_{pj}$. }
    \label{tab:fullperform}
\end{table}
One item of note is that TARES is the optimal opener in fewer cases now, with the mysterious LARES taking its place for most all of the positive $p$ values used. However, TARES still takes the prize for least expected rounds to win at $p=-0.5$, which is also one of the best performing $p$ values in terms of least losses.  It is interesting that, despite the dictionary in this case being much bigger and giving clearly different optimal openers for many $p$ values, it is still $p=-0.5$ that is the overall winner here.   Finally, apparently $\pm 2 \neq \pm\infty$ for this dictionary, as the asymptotic optimal openers are VENAL and SERAI as $p\to\pm\infty$, respectively.

\subsection{A more elaborate method}
While the $p$-optimality method above is very computationally fast and works seemingly pretty well for wise choices of $p$, there is at least one aspect of the method that is a bit lacking.  Namely, the method is only concerned with the possible lengths of the next viable set, and pays no mind to any other qualities that the next viable set might have. While it is clear that this is a weakness, it is worth exploring what kinds of qualities of potential viable sets might be of interest, to help guide the development of a more elaborate technique.

To begin this discussion, I would like to define what I have come to call a Fully Discernible Set (FDS).  A fully discernible set is any viable set such that there exists at least one ``key word'' (row in the matrix $C^m$) such that in that row, all of the words in the viable set will return a different color code.  I call this a fully discernible set because, if you were to guess the key word, the resulting color code would immediately allow you to discern what the correct answer is, fully (the single column in that row that contains that color code).  Using some of the terminology above, the current viable set is an FDS iff there is a row $i'$ of the matrix $C^m$ such that $M_{i'}(C^m)=N_v$, and in this case the word corresponding to row $i'$ is a key word to the FDS.  Given this definition, note that the maximum possible size for an FDS is $M$, as that is the largest possible value for $M_i(C^m)$.

Among all FDS, there are two general varieties: Internal-FDS (I-FDS) and External-FDS (E-FDS).  For a viable set to be an I-FDS simply means that there exists a key word to that FDS that is a member of the FDS.  The prototypical case: all viable sets of size two are I-FDS, with both members of the set as key words, as choosing either results in immediate victory or perfect knowledge of the correct answer.  Conversely, an E-FDS is an FDS where all the key words are not members of the FDS.  

When playing a game of Wordle, suppose you find that the current viable set is an I-FDS.  Again under the assumption that all words are equally likely to be the secret answer, then it is quite clear that the best possible thing you can do at that point is choose the (internal) key word.  This will cause you to win immediately with probability $1/N_v$, and if not you will win in the next round for sure.  But what about for an E-FDS?  Well, you generally have two plausible choices: pick the key word (which cannot cause you to win immediately, since it is external to the viable set) and then win in the next round for sure; or pick a member of the set instead and take your chance of winning right now, and even if you don't you will be guaranteed to still be in an FDS on the next round, since every subset of an FDS is also clearly an FDS.  This is perhaps most tempting if the size of the E-FDS is as small as possible (3), thus maximizing your chance of immediate win.  But, even in this case, in expectation this chancy choice is not better than simply choosing the external key word.  By choosing the external key word you will certainly win in 2 rounds.  By choosing an internal set member you win in one round with probability 1/3, otherwise you are left with an I-FDS (a set of size 2), meaning you win in two rounds with probability 1/3 and three rounds with probability 1/3, leading to an expected 2 rounds to win.  Any E-FDS with $N_v>3$ cannot be better than this, so you couldn't do better in expectation than choosing the external key in those cases either.  Hence, we will just make the choice that when an E-FDS is encountered, you take the sure bet and just choose the key word, in which case it is, regardless of its size, roughly on par with a viable set of length three.

One can of course extrapolate from the concept of an FDS to higher-orders.  For example, suppose my current viable set is not an FDS, but there does exists a row $i$ of $C^m$ such that all of the possible $\vec{v}^{m+1}$ sets that might arise should I choose $i$ are themselves FDS.  One might then call the current viable set an order 1 Discernible Set (1DS), which could be internal or external.  Similarly, one could conceive of an order $k$DS, such that it is not an FDS, but there exists a word such that all of the possible viable sets $\vec{v}^{k+1}$ arising from that word are themselves $(k-1)$ or lower order $DS$.  This is then starting to move away from a greedy algorithm, as we are inherently looking multiple steps ahead when discussing $k$DS.  For this reason, computing whether or not a set is $k$DS quickly becomes infeasible as $k$ grows from 1.  There is one freebie, though: a viable set of length $N_v=3$ is, if not itself an FDS, automatically an I-1DS.  This is because by picking any member of that set, you either get the answer correct immediately or you are left with a length 2 viable set, which is an I-FDS.

I will note that there is a nice relationship between the idea of an FDS and the $p$-optimal strategy.  Specifically, if a set is an FDS, then it is automatically the case that the key word $i$ will yield the $p$ optimal answer for any $p\neq 0$, since every $L(c,i,C^m)$ for that word will be 1, which will cause $\overline{f}(p)_i$ to be the smallest (or largest, depending on the sign of $p$) that it can possibly be.  Further, if the set is an I-FDS, then the internal key word will be the optimal choice, since we always check if there is a minimizer within the current viable set and pick that if so. 

Having defined a bunch of things, let me now describe how one might use them to try to obtain a modified optimal Wordle strategy.  Given the discussion above, in which we have noted that an E-FDS is essentially equivalent to a set of size three, and an I-FDS is akin to a set of size two (but maybe not quite as good, depending on its size), we could modify our $p$-optimal algorithm to take these facts into consideration.  Specifically, we first, as always, compute for each row $i$ of $C^m$ the lengths of all potential viable sets by color code returned, $L(c,i,C^m)$.   Then for each such color code that appears in row $i$ let us define a new value $\tilde{L}(c,i,C^m)$ with the following definition.  If $L<3$ or $L>M$ then $\tilde{L}=L$.  If $3\leq L(c,i,C^m)\leq M$, then the set is checked to determine if it is an FDS of some sort.  If the set is not an FDS then $\tilde{L}=L$.  If the set is an E-FDS then let $\tilde{L}=3$, to capture the fact that it is essentially the same as a length three set in terms of rounds to win.  Finally, if the set is an I-FDS then let $\tilde{L}=\frac{2}{L}2+\left(1-\frac{2}{L}\right)3$, to capture the fact that the expected number of rounds to win for this I-FDS is equal to this particular weighted average of the number of rounds needed to win for sets of size two and three.  Then we can simply modify our $p$-optimal algorithm to use $\tilde{L}^p$ rather than $L^p$, to obtain our $p$-FDS algorithm:
\begin{equation} \label{eq:pfdsopt}
\overline{\textrm{FDS}}(p)_i=\frac{1}{N_v-1_{M\in C^m_i}}\sum_{c=1,c\in C^m_i}^{M-1} L(c,i,C^m)\tilde{L}^p(c,i,C^m)~,
\end{equation}
where we again will minimize or maximize $\overline{\textrm{FDS}}(p)_i$ over $i$ depending on whether $p$ is positive or negative, respectively.

Now, all this is easy enough to state on paper, but starts to get a little bit slow, computationally.  Luckily, for the smaller dictionary at least, calculations using the $p$-FDS algorithm are just this side of bearable\footnote{In an early, quite slow version of my code, it was expected to take around nine days to compute all possible games of Wordle using $p$-FDS.  I set it to run in the background on my computer at work and moved on with my life.  Unfortunately, one week into the calculation my office building suffered a very rare power outage and all results were lost.  I interpreted this tragedy as a sign, and proceeded to optimize my code, which can now do the same calculation in roughly three hours.}, so I can provide some results.  For the value $p=-0.5$, chosen for its good performance noted above, the $p$-FDS optimal opening word is, you guessed it, TARES.  What a great word.  Further, I have simulated all possible Wordle games on normal mode for this specific $p$ value.  The mean number of rounds to win in this scenario is 3.7696, which is marginally better than the 3.7739 obtained from standard $p$-optimality for $p=-0.5$.  The $p$-FDS method does not have any losses at the seven round level, unlike the $p$-optimal method, which had two.  But, $p$-FDS does fail for six words at the 6.5 round level, as opposed to the $p$-optimal value of two words; but again note that failure at the 6.5 round level is only a possible failure, not a guaranteed one.  For those curious, these six possible fail words are BAKES, BASES, FAKES, HAKES, WAKES, and WASES, all following the \_A\_ES pattern noted above.  When comparing words directly, $p$-FDS wins in strictly fewer rounds than $p$-optimality for 742 words, while $p$-optimality wins in strictly fewer rounds than $p$-FDS for 702 words; the rest are ties.  So, in my mind it is clear that $p$-FDS takes the crown here, with the small caveat that there is a potential failure for four more words than $p$-optimality (but no guaranteed failures).


\section{So, we should always use TARES first, right?}
Maybe.  I think the answer to this depends on how similar your individual brain is to a computer.  Probably it is not so similar, so it will not be able to use all kinds of information equally well.  For instance, in my own personal play of Wordle, I find that the most helpful guesses are those that return color codes with G and Y.  This is because it is easier for me to start listing off words that I know contain certain letters than it is to start listing off words that \emph{don't} contain certain letters.  The algorithms above of course do not make any distinction here (aside from my tie breaking, which was chosen with this preference in mind), so they may suggest words that are not as useful to a human actually playing the game (assuming the human will not continue using a solver throughout).

Given this, it would not be unreasonable to choose an opening word that will contain, on average, the largest number of G and Y colors within the color code returned.  It is easy, given the precomputed $C^1$, to determine which word does this.  The answer is AROSE, with an average of 1.924 G plus Y.  The magical TARES is tied for number 4 with five other words by this accounting, with an average of 1.894 G plus Y.  This is pretty close to that of AROSE, so maybe TARES\footnote{Hey, I just noticed you can type TARES with only one hand.} is still the best, considering all the results above as well.

For the sake of some completeness, I will also mention here the words that have the highest average number of G colors, as well as highest average number of Y colors; these are SANES (with 0.9255 G on average) and RESAT (with 1.5712 Y on average).  On these lists, TARES appears at locations 19 and 849, respectively.  Wow, 849, huh?

Another thing to consider is some of our assumptions made in Section 1.  As already noted, Wordle most certainly does not choose a random secret answer from the full dictionary every day. I am also quite certain that the subdictionary that Wordle selects possible answer words from contains generally ``common'' words, so as to not upset players too much.  In my time playing Wordle thus far, probably the most ``obscure'' (but certainly not as obscure as PZAZZ) words I've seen as answers were ABBEY and KNOLL.  Of course, I would hope most all English speakers who are playing Wordle are familiar with these two words, even if they aren't the first to leap to mind when considering five letter words.  But there would maybe be a revolt if the answer for Wordle one day was YUKKY.  Yes, YUKKY.  As in ``Fozzie Bear's comedy is very YUKKY.''\footnote{At least, that's how I would use it.  Online dictionaries \cite{dict} seem to just list this as an alternate spelling for YUCKY, which I think is missing a golden linguistic opportunity.}

Of course, if at some point it becomes known that Wordle actually generated or generates its subdictionary of answers via some probability distribution $q_j$, where $q_j$ is the probability that dictionary word $j$ is chosen as the answer for a given game, the algorithms above can be adapted.  All one need do at each guess level $m$ is compute the current (or speculative, if you are doing this for a potential next viable set) probability $q_j^m$ for each word in the relevant viable set by simply using the known $q_j$ for each word, and then normalizing each so that they sum to one.  Hence
\[
q_j^m = \frac{q_j}{\sum_{j\in \vec{v}^m}q_j}
\]
Then our formulas above would use, in place of current probabilities $L(c,i,C^m)/N_v$, values
\begin{equation}
W(c,i,C^m)=\sum_{j,C_{ij}^m=c} q_j^m~.
\end{equation}
Of course other modifications might also make a lot of sense in this scenario.  For instance, two sets of the same size may not be very equivalent anymore if they have very different distributions of probabilities for the words they contain.  As an example, suppose one potential viable set of four words had a probability distribution that was very heavily weighted toward just a single word in the set, while another potential viable set of four words had equal probability for all four words.  In that case, you would certainly prefer the set with the unequal probability distribution, as that would make it more like a set of just one word than four.  So, in adapting the algorithms to these more general probabilities per word, likely a better choice than $L^p$ would make sense.  Unfortunately, such an exploration will have to wait until another day, or I will never complete this note.  This discussion naturally leads to our concluding section.

\section{Possible Extensions and Future Work}
All academic papers end with a section like this, so this one will too, despite the fact that it's not a normal academic paper.

One natural potential extension that follows directly from this work is to construct a time-varying strategy that is not constant in each round of the game.  That is, one can clearly (and efficiently) compute the  $p$-optimal word at each guess level for many different $p$ simultaneously, and then try to use some kind of heuristics to determine which of these to actually use for that guess.  I have not explored this possibility in enough detail at this time to say much intelligent about it, but it is something I think is an obvious next step.

Of greater difficulty, but also potentially greater payoff, would be to look into methods that are more ``global'' in nature than those above, and consider multiple steps as a time, rather than just one.  The $p$-FDS algorithm is a first small step in this direction, as it does effectively look ahead for a certain class of viable sets where the outcome is easily predicted in advance, but much more remains to be done here.  Again, having not thought about this in too much depth as yet, it may not be so horrible to look two steps ahead.  Given that, as the results above indicate, $p$-optimality seems to typically solve Wordle (technically my own Wordle implementation that definitely complies with the assumptions we have made in Section 1) in typically three or four steps, looking two steps ahead gets you a good way toward the solution, and might be quite useful.

Another extension would be to apply the methods here to alternate Wordles that use words of lengths other than five.  Of course, given the length chosen, the methods above may be completely infeasible from a computational point of view, as the number of possible words could be much bigger.  For example, given my source for dictionary words, there are 8,459 six letter words; 11,934 seven letter words; 13,057 eight letter words; 12,009 nine letter words; and the trend appears downward from there.  On the other end of the spectrum, there are 2,593 four letter words; 667 three letter words; and 77 two letter words (that seems high to me...).  These lower letter-count words could prove an interesting testing ground to see how well certain algorithms correctly find the true optimum, as it would be more plausible (but still maybe not so easy) to find the true global optimum by brute force on these much shorter lists (of fewer letters each).

Finally, consider again the possibility discussed above that the answers for any given Wordle be drawn from a given probability distribution. If the designer/maintainer of Wordle were so inclined, they could attempt to design this probability distribution in such a way as to optimally thwart my algorithmic solvers (or whatever other solvers might be out there).  This would pit the Wordle designers/maintainers and people such as myself in a full-on game-theoretic battle, which could be very interesting to attempt to solve, on both ends.





\bibliography{wordle.bib}{}
\bibliographystyle{IEEEtran}

\end{document}